\documentclass{amsart}
\usepackage{citesort}
\newtheorem{theorem}{Theorem}[section]
\newtheorem{lemma}[theorem]{Lemma}
\newtheorem{corollary}[theorem]{Corollary}
\theoremstyle{definition}

\theoremstyle{remark}
   \newtheorem{remark}[theorem]{Remark}
\numberwithin{equation}{section}
\newcommand{\IP}{\mathrm{P}}
\newcommand{\IE}{\mathrm{E}}
\begin{document}
\title{Brownian Sheet and Quasi-Sure Analysis}
\author{Davar Khoshnevisan}
\address{Department of Mathematics, University of Utah,
   Salt Lake City, UT 84112--0090}
   \thanks{Research supported in part by a grant from
   the NSF}
\email{davar@math.utah.edu}
\date{August 30, 2002}

\dedicatory{To Professor Miklos Cs\"org\H{o}
   on the occasion of his 70th birthday.}

\begin{abstract}
   We present a self-contained and
   modern survey of some
   existing quasi-sure results
   via the connection to the Brownian sheet.
   Among other things, we prove that
   quasi-every continuous function:
   (i) satisfies the local law of the iterated
   logarithm; (ii) has L\'evy's modulus of
   continuity for Brownian motion; (iii)
   is nowhere differentiable; and (iv)
   has a nontrivial quadratic variation. 
   We also present
   a hint of how to extend (iii) to
   obtain a quasi-sure refinement of the 
   M.\@ Cs\"org\H{o}--P.\@ R\'ev\'esz
   modulus of continuity for almost every continuous 
   function along the lines suggested by M.\@ Fukushima.
\end{abstract}

\maketitle

\section{Introduction}

Throughout, we let $\Omega$ denote the space of all continuous
functions $f:[0,1]\to\mathbb{R}$. As usual,
$\Omega$ is endowed with the compact-open topology
[i.e., the topology of uniform convergence] and its corresponding
Borel sigma-algebra $\mathcal{B}(\Omega)$. 
Then a number of classical, as well as
modern, theorems of probability theory can be interpreted
a saying something about the ``typical'' function in $\Omega$
in the following sense:\footnote{
The Baire category theorem provides us with another
notion of ``a typical function $f\in\Omega$,'' 
with a rich and colorful history dating back
to the works of Baire, Borel, and Lebesgue.
Among many gems of thought,
J.-P.\@ Kahane's article (\cite[\S3]{kahane}) contains a 
delightful discussion of the history of this subject.}
If we endow
$(\Omega,\mathcal{B}(\Omega))$ with the standard 
Wiener measure we then obtain the
\emph{classical Wiener space}, and various
``probabilistic'' results hold for almost 
every $f\in\Omega$. Here and throughout,
``almost every''
is tacitly understood to hold with respect 
to Wiener's measure
on $(\Omega,\mathcal{B}(\Omega))$.  As two notable
examples we can consider the following, although
frequently one 
thinks of these as statements about the Brownian motion:
\begin{enumerate}
   \item \label{lil}
      (Khintchine~\cite{khintchine}).
      \emph{Almost every $f\in\Omega$ 
      satisfies the local law of the iterated
      logarithm; i.e.,}
      \begin{equation}
         \limsup_{x\to 0} \frac{f(x)}{\sqrt{2x\ln
         \ln \left(\frac 1x\right)}}=-
         \liminf_{x\to 0} \frac{f(x)}{\sqrt{2x\ln
         \ln \left(\frac 1x\right)}}=1.
      \end{equation}
   \item \label{pwz}
      (Paley, Wiener, and Zygmund~\cite{paley-et-al}). 
      \emph{Almost every
      $f\in \Omega$ is nowhere-differentiable.}
\end{enumerate}
Other examples abound. 

One might ask for a 
more restrictive notion of what it means for
$f\in\Omega$ to be ``typical.'' In the
appendix to \cite{meyer}, D.\@ Williams
has proposed this problem, and has shown us an
interesting, less restrictive, class of 
``typical functions'' that is motivated
by infinite-dimensional diffusion theory,
and in particular the work of P.\@ Malliavin
in the said area (\cite{malliavin}).

Let $W$ denote the
two-parameter Brownian sheet, based on which
we can construct the 
\emph{Ornstein--Uhlenbeck Brownian sheet},
\begin{equation}\label{U}
   U(s,t) := e^{-s/2} W(e^s,t),
   \qquad\forall s,t\ge 0.
\end{equation}
One can think of this two-parameter process
as the ``evaluations'' of the
following infinite-dimensional 
(in fact, $\Omega$-valued) stochastic process that
is called the \emph{Ornstein--Uhlenbeck
process in Wiener space}:
\begin{equation}\label{Y}
   Y_s := U(s,\bullet),\qquad\forall s\ge 0.
\end{equation}
It is not difficult to see
that $Y$ is an $\Omega$-valued diffusion;
this follows at once from Theorem~\ref{th:SMP}
below, for instance. Moreover, since
$Y_0 =W(1,\bullet)$ is a standard Brownian motion,
it follows that the process $Y$ is a stationary 
diffusion on the space of continuous function,
and the invariant measure of $Y$ is Wiener's
measure.

Next, consider the hitting probabilities $\text{\rm Cap}(\bullet)$
of the diffusion $Y$ killed at an independent mean-one
exponential random variable; i.e., for any
Borel set $G\subset\Omega$, 
\begin{equation}
   \text{\rm Cap}(G):= \int_0^\infty
   e^{-s} \IP \left\{ 
   \exists s\ge 0:\ Y_s\in G\right\}
   \, ds.
\end{equation}

Following D.\@ Williams, we then say
that a Borel measurable set $G\subseteq\Omega$ holds
\emph{quasi-surely} if $\text{\rm Cap}(G^\mathsf{c})=0$.

It is not difficult to see that the set function
$\text{\rm Cap}$ is a natural capacity in the sense
of G.~Choquet. From this it follows that $G$ holds 
quasi-surely if and only
if its complement is a capacity-zero set; 
i.e., it is almost-surely never visited by
the Ornstein--Uhlenbeck process on $\Omega$.
Equivalently---and this requires only a 
moment of reflection---$G$ holds quasi-surely if and
only if 
\begin{equation}
   \IP\{ \forall s\ge 0:\ U(s,\bullet)\in G\}=1.
\end{equation}
Thanks to (\ref{U}), the quasi-sure analysis of
subsets of $\Omega$ can be related to the Brownian sheet.

An alternative, more direct, approach was proposed by
M.\@ Fukushima (\cite{fukushima}) who used the properties
of the Dirichlet form associated with the infinite-dimensional
process $Y$ to produce interesting quasi-sure theorems.
This was an exciting new development on the intersection of
probability and infinite-dimensional analysis,
and has led to a rich body of works;
cf.~\cite{CKS:00,CKS:99,CR,DM,DW1,DW2,%
fitzsimmons,fukushima,KD,KRS,KS,%
Komatsu-Takashima:84a,Komatsu-Takashima:84b,kono:84,Lyons:86,%
meyer:82,mountford:92,mountford:89,OP,penrose:90,penrose:89,%
shigekawa:84,takeda,walsh:86,walsh:84,zimmerman}. (Not all
of these references employ the quasi-sure notation 
in their presentation.)

The said connection to the Brownian sheet makes it
clear that whenever $G\subseteq\Omega$ holds quasi-surely,
then $G$ holds almost-surely as well. For a converse, it
has been noted in~\cite[p. 165]{fukushima}
that there are events that hold almost-surely and not
quasi-surely. For instance, consider $G$ to be
the collection of all continuous functions
$f:[0,1]\to\mathbb{R}$ such that $f(1)\neq 0$.
It is clear then that $G$ holds almost surely;
equivalently, with probability one a Brownian motion $B$
satisfies $B(1)\neq 0$. On the other hand, $G$ does
\emph{not} hold quasi-surely. [This is equivalent to the statement
that the Brownian sheet $W$ satisfies
$W(t,1)=0$ for some $t\ge 1$, which happens with
probability one
since $t\mapsto W(t,1)$ is a Brownian motion,
and hence is point-recurrent.]

I will say a few things in the final section
of this paper about the 
aforementioned analytical methods and their
potential-theoretic connections in turn.
However, this paper
is chiefly concerned with the aspects of quasi-sure analysis
that are close in spirit to what I believe may be
the general theme of this volume; namely, methods
that are based on finite-dimensional processes,
concentration, and
Gaussian inequalities.

On a few occasions, $\Omega$ will denote the
space of continuous functions $f:[0,1]\to\mathbb{R}^d$, and
$W$ will denote $d$-dimensional two-parameter
Brownian sheet, where $d\ge 1$. However, this should not cause
any confusion.

\section{A Strong Markov Property}\label{sec:SMP}
The following is an infinite-dimensional
strong Markov property of
the Brownian sheet. It is not a particularly
difficult result, but it is useful. 
In addition, this is a 
natural place to start our discussion.

Let $\mathfrak{F}:= \{ \mathfrak{F}(s);\, s\ge 0\}$ 
denote the filtration of $\sigma$-algebras defined
as follows:
For any $s\ge 0$, we first define
$\mathfrak{F}_{_{00}}(s)$ to be the $\sigma$-algebra
generated by the random variables
$\{W(r,t);\ r\in[0,s],\ t\ge 0\}$.
To each $\mathfrak{F}_{_{00}}(s)$
we can add all the $\IP$-null sets
and call the resulting $\sigma$-algebra
$\mathfrak{F}_{_0}(s)$. Finally, we make this
completed filtration right-continuous
in the usual way; namely, we let
$\mathfrak{F}(s):= \cap_{u>s}\mathfrak{F}_{_0}(u)$.

\begin{theorem}[A Strong Markov Property]
   \label{th:SMP}
   If $S$ is a finite $\mathfrak{F}$-stopping time,
   then the process 
   $t\mapsto W(S,t)$ is measurable
   with respect to $\mathfrak{F}(S)$. Moreover,
   the infinite-dimensional process $s\mapsto
   W(S+s,\bullet)-W(S,\bullet)$
   is totally independent of $\mathfrak{F}(S)$,
   and has the same law as $W$.
\end{theorem}

\begin{remark}
   This is a simple consequence of J.\@ B.\@ Walsh's
   strong Markov property with respect to weak
   stopping points; cf.~\cite[Theorem 3.6]{walsh:86} 
   or~\cite[Theorem 1.6]{walsh:84} for details.
\end{remark}

\begin{proof}
   Let $I(j;n)$ denote the half-open
   interval
   $[ j2^{-n},(j+1)2^{-n})$, and 
   for any fixed real number $r>0$, define 
   \begin{equation}
      S_{n,r}:= \sum_{j=0}^{\lfloor
      2^n r \rfloor} j2^{-n}
      \mathbf{1}_{ I(j;n) }(S).
   \end{equation}
   Since $S$ is an $\mathfrak{F}$-stopping time,
   so is $S_{n,r}$ for any fixed $n,r$;
   moreover, we have $S_{n,r}\le S\wedge r$, and
   as $n\uparrow+\infty$, then
   $S_{n,r}\uparrow S\mathbf{1}_{\{ S\le r\}}$.
   Now 
   \begin{equation}
      W(S_{n,r},t)=\sum_{j=0}^{\lfloor
      2^n r \rfloor} W(j2^{-n},t)
      \mathbf{1}_{ I(j;n) }(S).
   \end{equation}
   In particular, for any $t_1,\ldots,t_k\ge 0$,
   the vector $(W(S_{n,r},t_i))_{1\le i\le k}$
   is $\mathfrak{F} (S_{n,r})$-measurable,
   which is another way to say that $W(S_{n,r},\bullet)$
   is $\mathfrak{F} (S_{n,r})$-measurable. Since
   $\mathfrak{F} (S_{n,r})\subseteq\mathfrak{F}(S)$,
   this shows that $W(S_{n,r},\bullet)$ is 
   $\mathfrak{F}(S)$-measurable. Let $n\uparrow\infty$
   and $r\uparrow\infty$ (along rationals), and use
   the path continuity of $W$ to see that
   $W(S,\bullet)$ is $\mathfrak{F}(S)$-measurable,
   as asserted. Suppose $\Phi$ is the random function
   \begin{equation}
        \Phi(u):=
        \prod_{i=1}^m \phi_i \Big(
         W(u+s_i,t_i) - W(u,t_i)
         \Big),
   \end{equation}
   where $\phi_1,\ldots,\phi_m$ are 
   bounded continuous
   functions, and $t_1,\ldots,t_m\ge 0$. Then,
   for any bounded $\mathfrak{F}(S_{n,r})$-measurable
   random variable $\xi$,
   \begin{equation}
      \IE \left\{ \Phi ( S_{n,r} )  \cdot \xi \right\}
      \quad= \sum_{j=0}^{\lfloor 2^n r\rfloor}
      \IE  \Big\{ \Phi ( j2^{-n} )  \cdot \xi 
      \mathbf{1}_{
      I(n;j)}(S) \Big\}.
   \end{equation}
   The term ``$\xi$ times the indicator function'' is 
   $\mathfrak{F}(j2^{-n})$-measurable since
   $S$ is a stopping time. Therefore, the stationary
   independent-increments property of Brownian sheet
   implies that
   \begin{equation}
      \IE  \left\{ \Phi ( S_{n,r} )  \cdot \xi \right\}
      = \IE \{\Phi(0)\}\cdot
      \sum_{j=0}^{\lfloor 2^n r\rfloor}
      \IE  \left\{  \xi \mathbf{1}_{
      I(n;j)}(S) \right\}=\IE \{ \Phi(0)\}\cdot \IE \{\xi\}.
   \end{equation}
   This shows that a.s.,
   \begin{equation}
      \IE  \left\{\Phi(S_{n,r})\, \Big|\,\mathfrak{F}(S_{n,r})
      \right\}
      = \IE \{ \Phi(0)\},
   \end{equation}
   which is the desired
   strong Markov property in the case that
   $S\equiv S_{n,r}$. For the general case, 
   we let $n,r\uparrow\infty$ along rationals,
   and use the fact
   that $\mathfrak{F}(S_{n,r})\uparrow \mathfrak{F}(S)$,
   and that $\Phi(S_{n,r})\to \Phi(S)$ boundedly,
   together with H.\@ F\"ollmer's multiparameter
   version
   of Hunt's lemma (\cite[Lemma 2.3]{follmer}), to see that
   $\IE \{\Phi(S)\,|\,\mathfrak{F}(S)\}= \IE \{\Phi(0)\}$. This
   completes our proof.
\end{proof}

The above immediately yields a 0--1
law of the following form:

\begin{corollary}[A Zero-One Law]
   \label{co:0-1}
   If $S$ is a finite $\mathfrak{F}$-stopping time,
   then the following $\sigma$-algebra is trivial:
   \begin{equation}
      \mathfrak{A}(S) := \bigcap_{s\in\mathbb{Q}_+}
      \sigma \Big\{ W(S+s,\bullet)-W(S,\bullet)\Big\}.
   \end{equation}
\end{corollary}

\begin{proof}
   This follows from Theorem~\ref{th:SMP}
   and R.\@ M.\@ Blumenthal's $0$-$1$ law (cf.~\cite{blumenthal}).
   However, we include an argument that we will need
   later on, but do not wish to repeat.
   Consider the infinite-dimensional process
   $X(s;\delta):= W(S+s,\bullet)-W(S+\delta,\bullet)$,
   as $s$ varies over $(\delta,\infty)$ and $\delta>0$
   is a fixed number.
   Thanks to Theorem~\ref{th:SMP} applied to
   the $\mathfrak{F}$-stopping times of the form 
   $r+T$ (where $r$ is nonrandom),
   for any fixed $\delta>0$, all
   $s_1,s_2,\ldots,s_n> \delta$, and
   for all bounded continuous functionals 
   $\phi_1,\ldots,\phi_m$,
   \begin{equation}
      \IE \left.\left\{ \prod_{j=1}^m \phi_j 
      \left( X(s_j;\delta) \right)
      \, \right|\, \mathfrak{X}(\delta)
      \right\} = 
      \IE  \left\{ \prod_{j=1}^m \phi_j 
      \left( X(s_j;\delta) \right)\right\},\qquad\text{a.s.},
   \end{equation}
   where $\mathfrak{X}(\delta)$
   is the $\sigma$-algebra generated by
   $\{W(S+u,\bullet)-W(S,\bullet);\, u\in [0,\delta)\}$. Since
   $\cap_{\delta\in\mathbb{Q}_+}\mathfrak{X}(\delta)
   =\mathfrak{A}(S)$, we can let $\delta\downarrow 0$
   and apply Hunt's lemma 
   (cf. C.\@ Dellacherie and 
   P.-A.\@  Meyer~\cite[Chapter V, p. 25]{DM:80}) to deduce that
   \begin{equation}\label{eq:trivial}
         \IE \left.\left\{ \prod_{j=1}^m \phi_j 
         \left( X(s_j;0) \right)
         \, \right|\, \mathfrak{A}(S)
         \right\} = 
         \IE  \left\{ \prod_{j=1}^m \phi_j 
         \left( X(s_j;0) \right)\right\},\qquad\text{a.s.}
   \end{equation}
   The restriction $s_i>\delta$ has been removed
   since $\delta$ has been allowed to go to zero
   while keeping the $s_i$'s fixed. Thus,
   (\ref{eq:trivial}) holds for all bounded continuous
   functionals $\phi_1,\ldots,\phi_m$ and all 
   $s_1,\ldots,s_m>0$. A monotone class argument
   shows that $\mathfrak{A}(S)$ is independent of
   itself, and is trivial as a result.
\end{proof}

\section{A Law of the Iterated Logarithm}
\label{sec:zimmerman}
\begin{theorem}[G.~J.~Zimmerman~{\cite[Theorem 3]{zimmerman}}]
    \label{th:zimmerman}
    Quasi-every $f\in\Omega$
    satisfies the law of the iterated logarithm.
\end{theorem}

\begin{remark}
   For an analytic proof see
   \cite[Theorem 4]{fukushima}.
\end{remark}
 
Equivalently, Zimmerman's LIL 
states that with probability one,
\begin{equation}
    \limsup_{t \downarrow 0}\frac{ U(s,t) }{
    \sqrt{2 t \log \log \left( \frac 1t \right)}} =1,
    \qquad\forall
    s\in[0,1].
\end{equation}
Of course, the point is that the null set in question
does not depend upon $s\in[0,1]$ (or for that
matter upon $s\ge 0$ by scaling).
It is easy to see that the preceding equation can be translated
to the following statement about the Brownian
sheet: With probability one,
\begin{equation}\label{eq:zimmerman}
    \limsup_{t \downarrow 0}\frac{ W(s,t) }{
    \sqrt{2st \log \log \left( \frac 1t \right)}} =1,
    \qquad\forall
    s\in[1,e].
\end{equation}
In the next two subsections we will prove this particular
reformulation of Theorem \ref{th:zimmerman}. 
Before getting on with proofs, I would like
to mention---without proof---the following
theorem of \cite{mountford:92}.
Recall that a function $G:\mathbb{R}_+\to\mathbb{R}$ 
is an \emph{upper function}
for a function $g:\mathbb{R}_+\to\mathbb{R}$ 
if there exists $x_0>0$ 
such that for all $x\in[0,x_0]$, $g(x)\le G(x)$.

\begin{theorem}[T.\@ S.\@ Mountford~\cite{mountford:92}]
    \label{th:mountford}
    An increasing function 
    $t\mapsto\sqrt{t}\phi \left( \frac 1t\right)$ is in the upper
    class of quasi-every function in $\Omega$ if
    and only if
    \begin{equation}\label{eq:mountford}
        \int_4^\infty \phi^3(t) e^{-\phi^2(t)/2} \, 
        \frac{dt}{t}<+\infty.
    \end{equation}
\end{theorem}
On the other hand, the upper class for almost
all continuous paths has a different characterization
that is described in the 
classic paper \cite{erdos}; it states the following: 

\begin{theorem}[P.\@ Erd\H{o}s~\cite{erdos}]
    \label{th:erdos}
    An increasing function 
    $t\mapsto\sqrt{t}\phi\left( \frac 1t\right)$ is in the upper
    class of almost every function in $\Omega$ if
    and only if
    \begin{equation}\label{eq:erdos}
        \int_4^\infty \phi (t) e^{-\phi^2(t)/2} \, 
        \frac{dt}{t}<+\infty.
    \end{equation}
\end{theorem}

To illustrate, for any $\alpha>0$ define
\begin{equation}
    \phi_\alpha(t) := \sqrt{2 \log\log \left(
    \frac 1t\right)
    +\alpha \log\log\log \left(
    \frac 1t\right) },\qquad\forall t>4,
\end{equation}
and note that when $\alpha\in(3,5]$, 
$\phi_\alpha$ is an upper
function for almost all continuous paths, 
but it is \emph{not} an upper function 
for quasi-all of them.

We conclude this section by describing our
proof of Theorem \ref{th:zimmerman}.

\subsection{Upper Bound}\label{sec:zimmermanUB}
As in Khintchine's classical proof of the 
law of iterated logarithm, we begin by verifying
the following half of (\ref{eq:zimmerman}):
With probability one,
\begin{equation}\label{eq:zimmerman-UB}
    \limsup_{t \downarrow 0}
    \frac{ W(s,t) }{
    \sqrt{2st \log \log \left( \frac 1t \right)}} 
    \le 1,\qquad\forall
    s\in[1,e].
\end{equation}

To prove this, we will need an infinite-dimensional
reflection principle that we state in the following
abstract form.

\begin{lemma}[The Reflection Principle]\label{le:reflection}
    If $B$ is a continuous Brownian motion in a separable
    Banach space $\mathbb{B}$, and if $\mathcal{N}$ is any
    seminorm on $\mathbb{B}$ that is compatible with the
    topology of
    $\mathbb{B}$, then for all $T,\lambda>0$,
    \begin{equation}
        \IP \left\{ \sup_{t\in[0,T]}
        \mathcal{N}(B(t)) \ge \lambda \right\}
        \le 2 \IP \left\{ \mathcal{N}(B(T))
        \ge\lambda \right\}.
    \end{equation}
\end{lemma}

\begin{proof} 
    We follow the original ideas of
    D.\@ Andr\'e and P.\@ L\'evy that were developed 
    for 1-dimensional Brownian motion.
    
    Consider the stopping time 
    \begin{equation}
        \sigma :=\inf\left\{ t>0:\, 
        \mathcal{N}(B(t))\ge\lambda
        \right\}.
    \end{equation}
    Since $\mathcal{N}$ is compatible with the topology of
    $\mathbb{B}$, and since $W$ is continuous, 
    $\mathcal{N}(B(\sigma))=\lambda$ on $\{\omega:\,
    \sigma(\omega)<+\infty\}$. Now
    \begin{eqnarray}
        \lefteqn{\IP\left\{ \sup_{t\in[0,T]}
            \mathcal{N}(B(t))\ge\lambda \right\} =
            \IP\left\{ \mathcal{N}(B(T))
            \ge\lambda\right\}+
            \IP\left\{ \sigma<T ~,~
            \mathcal{N}(B(T))<\lambda\right\}}
            \nonumber\\
        & =& \IP\left\{ \mathcal{N}(B(T))
            \ge\lambda\right\}+
            \IP\left\{ \sigma<T ~,~
            \mathcal{N}(B(T)-B(\sigma)+B(\sigma))<
            \lambda\right\}\\
        & =& \IP\left\{ \mathcal{N}(B(T))
            \ge\lambda\right\}+
            \IP\left\{ \sigma<T ~,~
            \mathcal{N}(-B(T) +2B(\sigma))<\lambda\right\},
            \nonumber
    \end{eqnarray}
    thanks to the symmetry and independent-increments (i.e.,
    the strong Markov; cf. Theorem~\ref{th:SMP}) 
    properties of $W$. But the seminorm
    property of $\mathcal{N}$ insures us of its subadditivity.
    Thus, on $\{\sigma<+\infty\}$ we have
    \begin{equation}
        \mathcal{N}(-B(T)+2B(\sigma))\ge
        2\mathcal{N}(B(\sigma)) - \mathcal{N}(B(T))
        = 2\lambda - \mathcal{N}(B(T)).
    \end{equation}
    The previous two displays, used in conjunction,
    prove the result.
\end{proof}

\begin{proof}[Proof of (\ref{eq:zimmerman-UB})]
    Fix $c,\theta>1$, and consider the measurable events
    \begin{equation}
        \mathsf{F}_n := \left\{\omega:\
        \exists s\in[1,e],\,
        \sup_{0\le t\le \theta^{-n}} W(s,t) \ge
        \sqrt{2cs \theta^{-n}\log\log \theta^n}
        \right\}.
    \end{equation}
    We can rewrite $\mathsf{F}_n$ as follows.
    \begin{equation}
        \mathsf{F}_n  = \left\{\omega:\
        \sup_{0\le t\le \theta^{-n}}
        \mathcal{N}(B(t)) \ge \sqrt{2c \theta^{-n}
        \log\log \theta^n}  \right\},
    \end{equation}
    where for all $f\in\Omega$,
    $\mathcal{N}$ is the seminorm
    \begin{equation}
        \mathcal{N}(f)  = \sup_{1\le s\le e} \frac{
            f(s)}{\sqrt{s}},
    \end{equation}
    and $t\mapsto B(t)$ is the following Brownian motion
    in the Banach space $\mathbb{B}$ of continuous 
    functions on $[1,e]$ endowed with its compact-open
    topology:
    \begin{equation}
        B(t) (s)  = \frac{W(s,t)}{\sqrt{s}}.
    \end{equation}
    It follows readily that
    all of the assumptions of the reflection principle 
    are verified in the present context; cf. Lemma
    \ref{le:reflection}. Thus, the latter lemma implies
    that
    \begin{equation}\begin{split}
        \IP \left( \mathsf{F}_n \right)
            & \le 2 \IP\left\{ \sup_{1\le s\le e}
            \frac{W(s,\theta^{-n})}{\sqrt{s}} 
            \ge\sqrt{2c \theta^{-n}
            \log\log \theta^n} \right\}\\
        & = 2 \IP \left\{ \sup_{0\le s\le 1}
            O(s) \ge \sqrt{2c \log\log \theta^n} \right\},
    \end{split}\end{equation}
    where $O$ denotes a one-parameter Ornstein--Uhlenbeck
    process. C.\@ Borell's inequality (\cite{borell})
    shows that as $n\to\infty$, we have the estimate
    $\IP \left\{ \mathsf{F}_n \right\}
    \le n^{-c+o(1)}.$
    Since $c>1$, $n\mapsto \IP\{ \mathsf{F}_n\}$
    forms a summable
    sequence in $n$; thus, by the Borel--Cantelli
    lemma, with probability one,
    eventually $\mathsf{F}_n$ does not occur.
    Equivalently, with probability one,
    \begin{equation}
        \limsup_{n\to\infty} 
        \frac{\sup_{0\le t\le \theta^{-n}}
        W(s,t)}{\sqrt{2 s \theta^{-n}
        \log\log \theta^n}}\le \sqrt{c},
        \qquad\forall s\in[1,e].
    \end{equation}
    Since $c,\theta>1$  are arbitrary, monotonicity arguments
    yield (\ref{eq:zimmerman-UB}).
\end{proof}

\subsection{Lower Bound}
Theorem \ref{th:zimmerman} now follows once we show that
\begin{equation}
    \label{eq:zimmerman-LB}
    \limsup_{t\downarrow 0} 
    \frac{ W(s,t) }{
    \sqrt{2st \log \log \left( \frac 1t \right)}} 
    \ge 1,\qquad\forall
    s\in[1,e].
\end{equation}
\begin{proof}[Proof of (\ref{eq:zimmerman-LB})]
   Fix four constants $\varepsilon,c\in(0,1)$, $\tau>0$,
   and $\theta>1$,
   and consider the events
   \begin{equation}
       \mathsf{E}_n
       := \left\{ \omega:\
       \forall s\in[\tau, \tau(1+\varepsilon) ],\ 
       \frac{W(s,\theta^{-n})-
       W(s,\theta^{-n-1})}{ \sqrt{2cs
       (\theta^{-n}-\theta^{-n-1})
       \log\log \theta^n}}\ge 1 \right\}.
   \end{equation}
   Evidently, independently of $\tau>0$,
   \begin{equation}
      \IP ( \mathsf{E}_n )
      = \IP \left\{ \inf_{1\le s\le 1+\varepsilon} 
      \frac{B(s)}{\sqrt{s}} \ge 
      \sqrt{2c \log\log \theta^n}
      \right\},
   \end{equation}
   where $B$ is a Brownian motion. Trivially, for any
   given $c'\in(c,1)$, and as $n\to\infty$,
   \begin{equation}\begin{split}
      \IP  &\left( \mathsf{E}_n \right)
         \ge \IP\left\{ \inf_{%
         1\le s\le 1+\varepsilon} 
         B(s) \ge \sqrt{2c (1+\varepsilon)
         \log\log \theta^n}
         \right\}\hskip1in\\
      & \ge \IP \left\{ \sup_{%
         1\le s\le 1+\varepsilon}
         \left| B(s)-B(1) \right| \le 1
         \right\} \IP \left\{ B(1) \ge
         \sqrt{2c'(1+\varepsilon)
         \log\log \theta^n }\right\}  \\
      &= n^{-c'(1+\varepsilon) + o(1) }.
   \end{split}\label{eq:lower}\end{equation}
   Now if we also insist that $c(1+\varepsilon)<1$,
   then we can arrange things so that
   $c'(1+\varepsilon)<1$. In this case, the independence
   of $\mathsf{E}_1,\mathsf{E}_2,\ldots$, used in conjunction
   with (\ref{eq:lower}) and the Borel--Cantelli lemma,
   shows that infinitely many $\mathsf{E}_n$'s occur
   with probability one. Consequently, as long
   as $c(1+\varepsilon)<1$, then outside one null
   set, the following holds simultaneously for all
   $s\in[\tau,\tau(1+\varepsilon)]$:
   \begin{equation}\label{eq:19}\begin{split}
      &\limsup_{n\to\infty} \frac{W(s,\theta^{-n} )}{
         \sqrt{2cs(\theta^{-n}-\theta^{-n-1})
         \log\log \theta^n}}
         \hskip1.9in\\
      &\ \ge 1 - \limsup_{n\to\infty} \frac{
         |W(s,\theta^{-n-1} ) |}{%
         \sqrt{2cs(\theta^{-n}-\theta^{-n-1})\log\log
         \theta^n}}.\hskip1.5in 
   \end{split}\end{equation}
   The already-proven upper bound (cf.
   \S\ref{sec:zimmermanUB})
   implies that a.s., and simultaneously for all
   $s\in[\tau,\tau(1+\varepsilon)]$,
   \begin{equation}\label{eq:20}\begin{split}
      &\limsup_{n\to\infty} \frac{
         |W( s,\theta^{-n-1} ) |}{%
         \sqrt{2cs(\theta^{-n}-\theta^{-n-1})\log\log
         \theta^n}}\hskip1.9in\\
      &\ = \frac{1}{\sqrt{\theta-1}}
         \limsup_{n\to\infty} \frac{
         |W( s,\theta^{-n-1} ) |}{%
         \sqrt{2cs \theta^{-n-1} \log\log
         \theta^{n+1}}} \le \frac{1}{\sqrt{c(\theta-1)}}.
   \end{split}\end{equation}
   Thus, by (\ref{eq:19}) and (\ref{eq:20}),
   a.s., and simultaneously for all
   $s\in[\tau,\tau(1+\varepsilon)]$,
   \begin{equation}\label{eq:21}\begin{split}
      &\limsup_{t\to 0} \frac{W(s,t)}{
         \sqrt{2cs t\log\log \left( \frac 1t \right)}}
         \ge \limsup_{n\to\infty} \frac{W( s,\theta^{-n})}{
         \sqrt{2c s \theta^{-n} \log\log 
         \theta^n}}\hskip1in\\
      &\ \ge \sqrt{ 1-\frac{1}{\theta} }\cdot
         \limsup_{n\to\infty}
         \frac{W( s,\theta^{-n} ) }{
         \sqrt{2c s (\theta^{-n}-\theta^{-n-1}) 
         \log\log \theta^n}}
         \ge \sqrt{1 - \frac 1\theta} -
         \sqrt{\frac {1}{c\theta}}.
   \end{split}\end{equation}
   Since $\theta>1$ is arbitrary, we can let
   $\theta\uparrow+\infty$
   along a rational sequence to deduce that if
   $c(1+\varepsilon)>1$, then almost surely,
   \begin{equation}
      \limsup_{t\to 0} \frac{W(s,t)}{
      \sqrt{2cs t\log\log \left( \frac 1t \right)}} \ge 1,
      \qquad \forall s\in[\tau,\tau(1+\varepsilon)],\
      \forall\tau\in\mathbb{Q}_+.
   \end{equation}
   Let $c\uparrow (1+\varepsilon)^{-1}$ along a rational 
   sequence to see that
   \begin{equation}
      \limsup_{t\to 0} \frac{W(s,t)}{
      \sqrt{2 s t\log\log \left( \frac 1t \right)}} \ge \sqrt{
      \frac{1}{1+\varepsilon}},
      \qquad \forall s\in[\tau,\tau(1+\varepsilon)],\
      \forall\tau\in\mathbb{Q}_+.
   \end{equation}
   Equation (\ref{eq:zimmerman-LB}) follows from this
   readily.
\end{proof}

\section{J.\@ B.\@ Walsh's Proof of Theorem~\ref{th:zimmerman}}

The argument that was used to derive 
Theorem~\ref{th:zimmerman} (essentially
due to G.\@ J.\@ Zimmerman) is quite natural,
and has other uses in quasi-sure analysis
as we shall see in the next section.
I now wish to present a different
derivation of Theorem~\ref{th:zimmerman}---due to
J.\@ B.\@ Walsh---that is elegant and short.
It also has striking
consequences on the ``propagation of 
singularities'' along the Brownian
sheet.
The main ingredient of Walsh's proof is
the celebrated ``section theorem'' of
\cite{meyer} that,
borrowing from the words of M.\@ Sharpe,
``is one of the prime achievements
of [stochastic analysis].'' See
\cite[p. 388]{sharpe}.

\subsection{P.-A.\@  Meyer's Section Theorem}

In order to describe a version of Meyer's section
theorem that is suitable for our needs, we need
to recall a few notions from the general theory 
of processes.

Let $(\Omega,\mathfrak{G},\IP)$ denote a
filtered probability space, where the filtration
$\mathfrak{G}
:=(\mathfrak{G}_t)_{t\ge 0}$ satisfies the
``usual conditions'' of stochastic analysis,
i.e., $\mathfrak{G}_0$ contains all the
$\IP$-null sets, and $\mathfrak{G}_t=\cap_{r>t}\mathfrak{G}_r$.
A stochastic process $\{X_t\}_{t\ge 0}$
is said to be
\emph{optional} if: (i) For all $t\ge 0$, 
$X_t$ is $\mathfrak{G}_t$-measurable; and 
(ii) $t\mapsto X_t(\omega)$ is right-continuous
with left-limits for each $\omega\in\Omega$.
The \emph{optional}
$\sigma$-algebra $\mathfrak{O}$
is the smallest $\sigma$-algebra of
subsets of $[0,\infty)\times\Omega$ that renders
optional processes measurable; i.e.,
$\mathfrak{O}$ is the $\sigma$-algebra generated
by all sets of the form
$\{ (t,\omega)\in[0,\infty)\times\Omega:\
X_t(\omega)\in A \}$ where 
$A\subseteq\mathbb{R}$ is measurable
and $X$ is an optional process.
Finally, a stochastic set
$\Gamma \subseteq [0,\infty)\times\Omega$
is \emph{optional} if it is measurable
with respect to $\mathfrak{O}$.

\begin{theorem}[P.-A.\@  Meyer~{\cite[Chapter IV, pp.
   84--85]{meyer}}]
   \label{th:section}
   If $F$ is an optional set, then for
   every $\varepsilon>0$, there exists a stopping
   time $T_\varepsilon:\Omega\to[0,\infty)$ such that
   $\IP\{ T_\varepsilon<+\infty\} \ge \IP\{\Pi (F)
   \}-\varepsilon$,
   where $\Pi$ is the natural projection
   from $[0,\infty)\times\Omega$ onto $\Omega$.
\end{theorem}

\subsection{J.\@ B.\@ Walsh's Proof.}
For all $s\in[0,1]$ define
\begin{equation}
   L_s := \limsup_{t\to 0^+}\frac{W(s,t)}{
   \sqrt{2 t\log \log \left( \frac 1t\right)}}.
\end{equation}
The standard law of the iterated logarithm implies
that for each fixed $s>0$,
$L_s=\sqrt{s}$, a.s. In particular, thanks
to Fubini's theorem,
\begin{equation}\label{zero-leb}
   \text{Leb}\left\{ s>0:\
   L_s \neq \sqrt{s}\right\} =0,\qquad
   \text{a.s.},
\end{equation}
where $\text{Leb}$ denotes Lebesgue's
measure on $\mathbb{R}$.
Our goal is
to show that $\IP\{\forall s>0:\ L_s=\sqrt{s}\}=1$.
Suppose to the contrary
that $\IP\{\exists s>0:\ L_s\neq
\sqrt{s}\}>0$. 
We will use the section theorem to obtain
a contradiction. To do so, we need to meet
the conditions of Theorem~\ref{th:section}.
Let $\mathfrak{G}_t:= \mathfrak{F}(t)$ (the
filtration of  \S\ref{sec:SMP}), that
we recall satisfies the usual conditions.
Let 
\begin{equation}
   F:= \left\{ (t,\omega):\
   L_t(\omega) \neq \sqrt{t}
   \right\}.
\end{equation}
Since $s\mapsto W(s,\bullet)$ is continuous
(in the space of continuous functions on
$[0,1]$ endowed with compact-open topology),
and since $\mathfrak{F}$ is generated by the
latter process, the stochastic set $F$
is optional. By the section
theorem (Theorem~\ref{th:section}), 
there would then exist a finite 
$\mathfrak{F}$-stopping time $S$,
such that $\IP\{ L_S \neq \sqrt{S}\}>0$.
(Technical remark: The process $s\mapsto W(s,\bullet)$
is not real-valued. However, by considering
all processes of the form $\int W(s,t)\, \mu(dt)$
where $\mu$ is a linear combination of
point masses, we can see that
$\mathfrak{F}$ is generated by continuous
real-valued processes, so that the section theorem
can be applied as stated.)
Without loss of generality, we can assume
that there exists a $\delta>0$ (fixed and
nonrandom), such that
$\IP\{ L_S <\sqrt{S}-\delta \}>0$.
(If for some $\delta>0$,
$\IP\{ L_S >\sqrt{S} + \delta \}>0$, then a 
similar argument can be invoked to
get a contradiction.)

Thanks to the strong
Markov property (Theorem~\ref{th:SMP})
and the usual LIL, for any $s>0$, there exists 
a null set off of which,
\begin{equation}
   L_{S+s} \le L_S +
   \limsup_{t\to 0^+} \frac{W(S+s,t)-W(S,t)}{
   \sqrt{2t\log\log \left( \frac{1}{t} \right)}}
   =L_S + \sqrt{s}.
\end{equation}
Therefore, for all $s>0$ sufficiently small, 
\begin{equation}
   \IP \left\{L_{S+s} < \sqrt{S+s}
   \right\} \ge \IP  \{ L_{S+s} <\sqrt{S} + \sqrt{s}
   -\delta \}>0.
\end{equation} 
We can integrate this $[ds]$
and use Fubini's theorem to deduce
that  with positive
probability,
\begin{equation}
   \text{Leb}
   \left\{  s>0:\ L_{S+s} <\sqrt{S+s} \right\}>0,
\end{equation}
which contradicts (\ref{zero-leb}).\hfill$\square$

\section{Modulus of Continuity}

A well-known result of P.\@ L\'evy (\cite{levy}) states that
almost all $f\in\Omega$ have
the following uniform modulus of continuity:
\begin{equation}\label{eq:modulus}
   \limsup_{\varepsilon\to 0^+} \sup_{%
   \scriptstyle u,v\in[0,1]:
   \atop\scriptstyle |u-v|\le\varepsilon}
   \frac{ \left| f(u) - f(v) \right|}{%
   \sqrt{2\varepsilon\log\left( \frac{1}{\varepsilon}
   \right)}}=1.
\end{equation}
In an elegant paper that popularized the
subject of quasi-sure analysis, M.\@ Fukushmia
proved the following quasi-sure analogue.

\begin{theorem}[M.\@ Fukushima~{\cite[Theorem 3]{fukushima}}]%
   \label{thm:fukushima}
   Quasi-every $f\in\Omega$
   has the uniform modulus of continuity
   described by (\ref{eq:modulus}).
\end{theorem}

The argument of \cite{fukushima} involves
infinite-dimensional analysis and Dirichlet form
estimates. Instead of going that route, we follow a more
classical route that has the advantage of providing
us with a more delicate result. To explain this
extension, we first recall that
in their book~(\cite{CR1}), M.\@ Cs\"org\H{o} and 
P.\@ R\'ev\'esz
have shown us that  even if we replace the
$\limsup$ by a proper limit there,
(\ref{eq:modulus}) holds for almost every function.
By adapting their argument, we plan to prove the
following refinement of Theorem 
\ref{thm:fukushima}.

\begin{theorem}[M.\@ Fukushima~{\cite[Theorem
   3]{fukushima}}]\label{thm:fukushima'}
   Quasi-every $f\in\Omega$ satisfies
   (\ref{eq:modulus}) with $\limsup$ replaced
   by a proper limit.
\end{theorem}

\begin{proof}
   Our goal is to show that with probability one,
   \begin{equation}
      \lim_{\varepsilon\to 0^+} \sup_{%
      \scriptstyle u,v\in[0,1]:
      \atop\scriptstyle |u-v|\le\varepsilon}
      \frac{ \left| U(s,u) - U(s,v) \right|}{%
      \sqrt{2\varepsilon
      | \log \varepsilon |}}=1,\qquad\forall s\in[0,1].
   \end{equation}
   Instead, we will prove the following stronger
   result: Almost surely, as $\varepsilon\to 0^+$,
   \begin{equation}\label{unif:goal}
      \sup_{%
      \scriptstyle u,v\in[0,1]:
      \atop\scriptstyle |u-v|\le\varepsilon}
      \frac{ \left| U(s,u) - U(s,v) \right|}{%
      \sqrt{2\varepsilon
    | \log \varepsilon |}} \longrightarrow 1,
   \end{equation}
   \emph{uniformly} for all $s\in[0,1]$.
   Clearly, this is equivalent to the following statement
   about the Brownian sheet that we propose to derive: 
   With probability one,
   \begin{equation}\label{unif:goal'}
      \lim_{\varepsilon\to 0^+} 
      \sup_{s\in[1,e]} \left| \sup_{%
      \scriptstyle u,v\in[0,1]:
      \atop\scriptstyle |u-v|\le\varepsilon}
      \frac{ \left| W(s,u) - W(s,v) \right|}{%
      \sqrt{2s\varepsilon
    | \log \varepsilon |}} -1 \right| = 0.
   \end{equation}

   \subsection{The Upper Bound}
   Fix $0<\theta<1$, and
   define $\delta_n :=n^2\theta^n$,
   and $\Theta_n := \{ j\theta^n:\,%
   0\le j\le\theta^{-n}\}$, and notice that
   \begin{equation}
      \IP\left\{ \max_{\scriptstyle
      u,v\in\Theta_n:\atop\scriptstyle
      |u-v|\le\delta_n}
      \left| W(s,u) - W(s,v) \right| \ge
      \sqrt{2s \delta_n \lambda} \right\} \le
      2 n^2 |\Theta_n| e^{-\lambda},
   \end{equation}
   since for every  $v\in\Theta_n$, there are
   no more than $2n^2$ many
   $u\in\Theta_n$ such that $|u-v|\le\delta_n$.
   Let $\lambda := p\log(\delta_n^{-1})$ for a fixed
   $p>1$, and appeal to our abstract form
   of reflection principle 
   (Lemma \ref{le:reflection}) to deduce that for
   any $\tau>0$,
   \begin{equation}\begin{split}
      &\IP\left\{ \max_{s\in[0,\tau]}
         \max_{\scriptstyle
         u,v\in\Theta_n:\atop\scriptstyle
         |u-v|\le\delta_n} 
         \left| W(s,u) - W(s,v) \right| \ge
         \sqrt{2 p \tau
         \delta_n \log(\delta_n^{-1})} \right\}\\
      &\quad\le 2 \IP\left\{
         \max_{\scriptstyle
         u,v\in\Theta_n:\atop\scriptstyle
         |u-v|\le\delta_n} 
         \left| W(\tau,u) - W(\tau,v) \right| \ge
         \sqrt{2 \tau p \delta_n \log(\delta_n^{-1})} 
         \right\}\\
      &\quad \le \theta^{n(p-1)+o(n)}.
   \end{split}\end{equation}
   On the other hand, by Kolmogorov's continuity theorem
   (cf. \cite[Chapter 5, Exercise~2.5.1]{MPP} 
   for a suitable version 
   of the latter theorem), for any integer $k\ge 1$,
   there exists a contant $A_k$ such that for all
   $r\in(0,1)$,
   \begin{equation}\label{eq:kolmogorov}
      \IE \left[ \sup_{s\in[0,e]}
      \sup_{\scriptstyle u,v\in[0,1]:\atop\scriptstyle
      |u-v|\le r} \left| W(s,u)-W(s,v) \right|^k \right]
      \le A_k r^{\frac k2}.
   \end{equation}
   Thus, thanks to the triangle inequality,
   for any fixed $\eta>0$ and $\ell\ge 0$,
   \begin{equation}
      \IP\left\{
      \sup \left| W(s,u)-W(s,v) \right| \ge \psi_n
      \right\} \le \theta^{n(p-1)+o(n)} + 2A_4 n^{-2},
   \end{equation}
   where the supremum is taken over all 
   $s\in [1+\ell\eta , 1+(\ell+1)\eta]$
   and $u,v\in[0,1]$ such that $|u-v|\le\delta_n$,
   and $\psi_n:=\sqrt{2p\{ 1+(\ell+1)\eta \}\delta_n\log 
   \left( \delta_n^{-1} \right)}+2n\theta^{\frac n2}$.
   The displayed probability is summable in $n$. But
   as $n\to\infty$, we have
   $\delta_n=(1+o(1))\delta_{n+1}$; moreover $p>1$
   is arbitrary. Therefore,
   the Borel--Cantelli lemma and monotonicity
   together show that with probability one
   \begin{equation}
      \limsup_{\varepsilon\to 0^+}\ \sup 
      \frac{ \left| W(s,u)-W(s,v)\right| }{%
      \sqrt{2\varepsilon 
      | \log \varepsilon |}}
      \le \sqrt{  1+(\ell+1)\eta  },
   \end{equation}
   where the supremum is taken over all
   $s\in [1+\ell\eta , 1+(\ell+1)\eta]$
   and $u,v\in[0,1]$ such that $|u-v|\le\varepsilon$.
   For the $s$ in question, $1+(\ell+1)\eta \le s(1+\eta)$.
   Since there are only finitely many integers $\ell$
   to consider (namely, $0\le\ell\le e/\eta$),
   \begin{equation}
      \limsup_{\varepsilon\to 0^+}\sup_{s\in[1,e]}
      \sup_{\scriptstyle u,v\in[0,1]:\atop\scriptstyle
      |u-v|\le\varepsilon}
      \frac{ \left| W(s,u)-W(s,v)\right| }{%
      \sqrt{2s\varepsilon | \log \varepsilon |}}
      \le \sqrt{ 1+\eta}.
   \end{equation}
   Let $\eta\downarrow 0$ along a rational sequence
   to deduce that a.s.,
   \begin{equation}\label{unif:UB}
      \limsup_{\varepsilon\to 0^+}\sup_{s\in[1,e]}
      \sup_{\scriptstyle u,v\in[0,1]:\atop\scriptstyle
      |u-v|\le\varepsilon}
      \frac{ \left| W(s,u)-W(s,v)\right| }{%
      \sqrt{2s\varepsilon | \log \varepsilon |}}
      \le 1.
   \end{equation}
   This proves half of (\ref{unif:goal'}).

   \subsection{The Lower Bound}
   Notice that for any $p\in(0,1)$ fixed, all
   integers $n\ge 1$, and all $s>0$,
   \begin{equation}\begin{split}
      &\IP\left\{ \sup_{\scriptstyle
         u,v\in[0,1]:\atop\scriptstyle |u-v|\le n^{-1}}
         \left| W(s,u) - W(s,v) \right| \le \sqrt{
         \frac {2ps}{n} \log n}
         \right\}\\
      & \quad \le \IP\left\{ \max_{%
         0\le j\le n-1}
         \left| W \left( s, \frac{j+1}{n}
         \right) - W \left( s, \frac{j}{n}
         \right) \right| \le \sqrt{
         \frac {2ps}{n} \log n}
         \right\}\qquad\qquad\ \ \\
      & \quad = \left( 1 - \IP \left\{
         |\mathcal{N}|  > \sqrt{2p\log n}\right\}\right)^n,
   \end{split}\end{equation}
   where $\mathcal{N}$ denotes a standard normal variable.
   Consequently, as $n\to\infty$,
   \begin{equation}
   \IP\left\{ \sup_{\scriptstyle
      u,v\in[0,1]:\atop\scriptstyle |u-v|\le n^{-1}}
      \left| W(s,u) - W(s,v) \right| \le \sqrt{
      \frac {2ps}{n} \log n}
      \right\} \le \exp\left(
      -n^{1-p+o(1)}\right),
   \end{equation}
   where $o(1)$ goes to $0$ uniformly in $s>0$.
   Let 
   \begin{equation}
      \sigma_n:= \left\{ 1+jn^{-4}:\
      0\le j\le (n+1)^4\right\}
   \end{equation}
   to see that as $n\to\infty$,
   \begin{equation}
      \sum_{n=1}^\infty
      \IP\left\{
      \min_{s\in\sigma_n}
      \sup_{\scriptstyle
      u,v\in[0,1]:\atop\scriptstyle |u-v|\le n^{-1}}
      \frac{ \left| W(s,u) - W(s,v) \right| 
      }{\sqrt{s}} \le \sqrt{
      \frac {2p}{n} \log n}
      \right\} <+\infty.
   \end{equation}
   Since $p\in(0,1)$ is arbitrary,
   by the Borel--Cantelli, and
   after  applying another interpolation
   argument involving (\ref{eq:kolmogorov})
   and yet another monotonicity argument,
   we conclude that a.s.,
   \begin{equation}
      \liminf_{\varepsilon\to 0^+} \inf_{s\in[1,e]}
      \sup_{\scriptstyle u,v\in[0,1]:\atop
      \scriptstyle |u-v|\le \varepsilon}
      \frac{\left|W(s,u)-W(s,v)\right|}{
      \sqrt{2s \varepsilon| \log \varepsilon |}} 
      \ge 1.
   \end{equation}
   Together with (\ref{unif:UB}), this yields
   (\ref{unif:goal'}) whence (\ref{unif:goal}).
\end{proof}

\section{Nowhere-Differentiability}

A classical result of R.\@ E.\@ A.\@ C.\@ Paley,
N.\@ Wiener and A.\@ Zygmund (\cite{paley-et-al})
states that almost all continuous functions
are nowhere-differentiable;
see also~\cite{dvoretzky-et-al}. This has been extended
in various directions in \cite{CR},
a consequence of which is the following;
see~\cite[Theorem 2]{fukushima} for 
an analytical proof of most of this theorem.

\begin{theorem}[M.\@ Cs\"org\H{o} and 
    P.\@ R\'ev\'esz~{\cite[Theorem 3]{CR}}]\label{th:CR}
    Quasi-every $f\in\Omega$ is 
    nowhere-differentiable.
\end{theorem}

We mention---without
proof---the following uniform modulus
of nondifferentiability that is a two-parameter
extension of the result of \cite{CR:79}: 
With probability one, as $\varepsilon\to 0^+$,
\begin{equation}
    \inf_{t\in[0,T]}
    \sup_{u\in[0,\varepsilon]}
    \frac{\left| U(s,t+u)-U(s,t)\right|}{%
    \sqrt{\varepsilon / |\log\varepsilon|}}
    \longrightarrow \frac{\pi}{\sqrt{8}}  ,
\end{equation}
\emph{uniformly} for all $s\in[0,1]$.
In particular, the Cs\"org\H{o}--R\'ev\'esz modulus of
nondifferentiability holds for quasi-all 
continuous functions. 

In fact, M.\@ Cs\"org\H{o} and 
P.\@ R\'ev\'esz~\cite[Theorem 3]{CR}
proved the much stronger theorem that the Brownian
sheet is nowhere-differentiable along \emph{any line 
in the plane.} Here we have specialized this
result to the simpler case where the lines 
are parallel to one of the axes. A consequence of this
more general theorem of~\cite{CR} is that the
level sets of the Brownian sheet a.s. do not contain
straight-line segments. More recently, 
R.\@ C.\@ Dalang and T.\@ S.\@ Mountford 
have discovered a striking generalization of this
fact: 
   
\begin{theorem}[%
   R.\@ C.\@ Dalang and T.\@ S.\@ Mountford~\cite{DM}]
   With probability one,
   the level curves of the Brownian sheet 
   do not contain
   any curve that is differentiable somewhere.
\end{theorem}

\begin{proof}[Proof of Theorem \ref{th:CR}]
    Motivated by the analysis of
    \cite{CR}, our strategy
    will be to show that
    if $T>0$ is held fixed, then with probability one,
    \begin{equation}\label{eq:4.1}
        \inf_{s\in[0,1]}\ \limsup_{n\to\infty}
        \inf_{t\in[0,T]}\
        \sup_{u\in[0,n^{-1}]} \frac{
        \left| U(s,t+u) - U(s,t)
        \right| }{n^{-1}} =+\infty.
    \end{equation}
    This implies Theorem \ref{th:CR}. With this in mind,
    let us first recall the following estimate: If $B$ denotes
    standard Brownian motion, then 
    there exists a constant
    $p>0$ such that for all $x>0$,
    $-\log \IP \left\{ \sup_{0\le t \le 1}
    |B(t)| \le x \right\}$ is bounded below by 
    $(p x^2)^{-1}$; cf. K.\@ L.\@ Chung
    \cite[Theorem 2]{chung}.
    Equivalently, for any $n\ge 1$ and all $t,a>0$,
    \begin{equation}
        \IP \left\{ \sup_{u\in[0,n^{-1}]}
        \left| B  \left( t+u \right) - B(t) \right| 
        \le \frac{a}{n} \right\} \le 
        \exp\left( -\frac{n}{p a^2} \right),
    \end{equation}
    for the same constant $p$ whose
    valued does not depend on our choice of $(a,t,n)$.
    [In fact, the optimal choice is $p= 8\pi^{-2}$.]
    Consequently,
    \begin{equation}\begin{split}
        &\IP \left\{ \inf_{t\in[0,T]}\
            \sup_{u\in[0,n^{-1}]}
            \left| B  \left( t+u \right) - B(t) \right| 
            \le \frac{a}{n} \right\} \\
        &\quad \le \IP \left\{ 
            \min_{0\le j\le n^4}\
            \sup_{u\in[0,n^{-1}]}
            \left| B  \left(  j n^{-4}T +u \right) - 
            B \left( jn^{-4}T \right)
            \right| \le \frac{2a}{n} \right\}\qquad\qquad\\
        &\qquad +\IP \left\{
            \sup \left| B(u) - B(v) \right|\ge
            \frac{a}{2n} \right\} \\
        &\quad \le (1+n^4) 
            \exp\left( -\frac{n}{4p a^2}
            \right) + \left( \frac{n}{2a}\right)^2
            \IE \left[ \sup
            \left| B(u) - B(v) \right|^2 \right],
            \qquad\qquad
    \end{split}\end{equation}
    where the last two suprema 
    are over all $u,v\in[0,1]$ such that
    $|u-v|\le n^{-4}T$.
    By (\ref{eq:kolmogorov}), the last expectation is
    seen to be no more than $A_2Tn^{-4}$.
    Consequently,
    \begin{equation}
        \sum_{n=1}^\infty
        \IP \left\{ \inf_{t\in[0,T]}\
        \sup_{r\in[0,n^{-1}]}
        \left| B  \left( t+r \right) - B(t) \right| 
        \le \frac{a}{n} \right\}  < +\infty.
    \end{equation}
    An interpolation argument improves this condition
    to the following one for the
    Brownian sheet $W$:
    \begin{equation}
        \sum_{n=1}^\infty 
        \IP \left\{ \inf_{s\in[0,e]}\
        \inf_{t\in[0,T]}\
        \sup_{r\in[ 0,n^{-1}]}
        \left| W  \left( s,t+r \right) - W(s,t) \right| 
        \le \frac{a}{n} \right\}<+\infty,
    \end{equation}
    from which we can easily deduce (\ref{eq:4.1}).
\end{proof}

\section{Quadratic Variation}

We now come to the theorem that started much
of the interest in quasi-sure analysis. Namely,
D.\@ Williams's quasi-sure refinement of the 
classical theorem of P.\@ L\'evy that
states that for almost every continuous function $f$, 
at time $t$ the function $f$ has finite quadratic 
variation $t$; cf. the appendix of~\cite{meyer:82}.

\begin{theorem}[D.\@ Williams {\cite[Appendix]{meyer}}]
   For each $n=1,2,\ldots$,
   let $0=\pi_{0,n}<\pi_{1,n}<\cdots<\pi_{n,n}=1$ denote
   a partition of $[0,1]$ such that 
   $\sum_n\max_{1\le j\le n}(\pi_{j,n}-\pi_{j-1,n})<+\infty$.
   Then quasi-every $f\in\Omega$ has
   the following property: For all $t\in[0,1]$,
   \begin{equation}
      \lim_{n\to\infty} \sum_{j=1}^n \Big|
      f(\pi_{j,n} t) - f(\pi_{j-1,n}t) \Big|^2 
      =t.
   \end{equation}
\end{theorem}

\begin{proof}
   In the interest of saving space, I will prove the slightly weaker
   statement
   that for each fixed $t>0$, with probability one,
   \begin{equation}\label{QV}
      \lim_{n\to\infty} \sup_{1\le s\le e}
      \left| \sum_{j=1}^n \big| W(s,\pi_{j,n} t)-
      W(s,\pi_{j-1,n} t) \big|^2 - st \right|=0.
   \end{equation}
   The proof that is to follow can be enhanced,
   using similar ideas, to show that in fact
   the above holds outside a single null set,
   uniformly for all $t\in[0,1]$, which yields the
   full statement of the theorem.
   
   We define the following, 
   all the time keeping $t\in[0,1]$
   fixed:
   \begin{equation}
      \theta_j(s)  := \Big| W(s,\pi_{j,n} t)-
         W(s,\pi_{j-1,n} t) \Big|^2.
   \end{equation}
   A few lines of calculations then show the existence of
   a universal constant $K_1$ such that
   $\IE \{ |\theta_j(s')-\theta_j(s) |^8 \}
   \le K_1 |s-s'|^4 (\pi_{j,n}-\pi_{j-1,n})^4.$
   Thus, by the Kolmogorov continuity theorem
   (\cite[Chapter 5, Exercise~2.5.1]{MPP}), we can find a
   universal constant $K_2$ such that for all $\eta\in(0,1)$,
   $j=1,\ldots,n$, and $n=1,2,\ldots$,
   \begin{equation}\label{eq:mod}
       \IE\left\{ \max_{\scriptstyle
       s,s'\in[1,e]:\atop\scriptstyle
       |s-s'|\le\eta}
       \left| \theta_j(s)-\theta_j(s') \right|^8 \right\}
       \le K_2 \eta^3 (\pi_{j,n}-\pi_{j-1,n})^4.
   \end{equation}
   Now choose an equipartition $\mathbb{S}_n$ of $[1,e]$ with
   $\text{mesh}(\mathbb{S}_n)\to 0$ at a rate to be
   described shortly, and note that for any
   fixed $\delta>0$,
   \begin{equation}
      \IP\left\{ \max_{s\in \mathbb{S}_n} \left| V_n(s,t)
      \right| \ge \delta \right\} \le \frac{\#(\mathbb{S}_n) 
      }{ \delta^2} \max_{s\in \mathbb{S}_n} \text{Var}( V_n(s,t)),
   \end{equation}
   where
   \begin{equation}
      V_n(s,t) := \sum_{j=1}^n \left|
      W(s,\pi_{j,n}t)-W(s,\pi_{j-1,n}t)\right|^2-st.
   \end{equation}
   On the other hand, $V_n(s,t)$ is a sum of $n$ i.i.d.
   random variables, and a simple computation yields
   a universal constant $K_3$ such that uniformly 
   for all $s\le e$,
   $\text{\rm Var}(V_n(s,t))\le K_3\sum_{j=1}^n
   (\pi_{j,n}-\pi_{j-1,n})^2\le K_3 \max_{1\le j\le n}
   (\pi_{j,n}-\pi_{j-1,n})):= K_3 \|\Pi\|_n.$
   This yields
   \begin{equation}\label{eq:max1}
      \IP\left\{ \max_{s\in \mathbb{S}_n} \left| V_n(s,t)
      \right| \ge \frac{\delta}{2} \right\} \le 
      \frac{K_3}{\delta^2}
      \#(\mathbb{S}_n) \|\Pi\|_n.
   \end{equation}
   Thus, for all $n$ large,
   \begin{equation}\label{eq:max2}\begin{split}
       \IP&\left\{ \sup_{s\in[1,e]} \left| V_n(s,t)
          \right| \ge \frac{\delta}{2} \right\} \cr
       &\le \frac{K_3}{\delta^2}
          \#(\mathbb{S}_n) \|\Pi\|_n
       + \IP \left\{ \max_{\scriptstyle s,s'\in
         [1,e]\atop\scriptstyle
         |s-s'|\le\text{mesh}(\mathbb{S}_n)} \sum_{j=1}^n
         \left| \theta_j(s)-\theta_j(s')\right|\ge \frac{\delta}{4}
         \right\}\cr
       &\le \frac{K_3}{\delta^2}
          \#(\mathbb{S}_n) \|\Pi\|_n 
          + K_2 \left(\frac{4}{\delta}\right)^8
          \left[ \text{mesh}(\mathbb{S}_n)\right]^3 \sum_{j=1}^n
          (\pi_{j,n}-\pi_{j-1,n})^4\\
       &\le \frac{K_3}{\delta^2}
          \#(\mathbb{S}_n) \|\Pi\|_n 
          + K_2 \left(\frac{4}{\delta}\right)^8
          \left[ \text{mesh}(\mathbb{S}_n)\right]^3 \|
          \Pi\|_n^3, \hskip1.3in
   \end{split}\end{equation}
   thanks to (\ref{eq:max1}).
   Since $\sum_n \|\Pi\|_n<\infty$
   we can always choose $\mathbb{S}_n$ such that
   the left-hand side of the preceding display is summable
   (in $n$), and this proves (\ref{QV}).
\end{proof}

\section{W.\@ S.\@ Kendall's Theorem}

Upto this point, we have adopted the viewpoint
that the Brownian sheet (equivalently,
the Ornstein--Uhlenbeck process) is a natural
diffusion on the space of continuous functions.
While this viewpoint provides us with a great
deal of insight about the sheet, it completely ignores
the effect of the geometry of the parameter space
on the process. The following is a delightful example of
the subtle effect of the geometry of the parameter
space, and was discovered by W.\@ S.\@ Kendall.

\begin{theorem}[W.\@ S.\@ Kendall%
   ~{\cite[Theorem 1.1]{kendall}}]\label{th:kendall}
   If $s,t>0$ are arbitrary but fixed,
   then with probability one, the
   level curve of $W$ that goes through $(s,t)$ is
   totally disconnected.
\end{theorem}

\begin{remark}
   One can restate Theorem~\ref{th:kendall} 
   in the following manner:
   Consider the level curve of $W$ that goes
   through $(s,t)$, and let
   \begin{equation}\label{CC}\begin{split}
      \Gamma(s,t) := &\text{ The connected
         component of }\Big\{ (u,v):\
         W(u,v)=W(s,t) \Big\}\\
      &\text{ that contains }(s,t).
   \end{split}\end{equation}
   By the continuity of $W$, this definition is
   perfectly well-defined, and 
   Theorem~\ref{th:kendall} asserts 
   that with probability one,
   $\Gamma(s,t)=\{(s,t)\}$. In other words,
   the excursion at the level $W(s,t)$ corresponding
   to the time-point $(s,t)$ is trivial. See
   the works of R.\@ C.\@ Dalang and J.\@ B.\@ Walsh~\cite{DW1,DW2}
   who present very precise descriptions of the local structure
   of the excursions of the Brownian sheet.
\end{remark}

\begin{remark}
   The null set in question depends on the
   choice of $(s,t)$. Moreover, a little thought shows
   that Theorem~\ref{th:kendall} cannot hold 
   a.s. simultaneously
   for all $(s,t)$ in any given open set. In this
   sense, this result is optimal.
\end{remark}

The proof of Kendall's theorem rests on the
following zero-one law; it is a two-parameter 
analogue
of the infinite-dimensional zero-law described
earlier in Corollary~\ref{co:0-1}.

\begin{lemma}[S.\@ Orey and W.\@ E.\@ Pruitt~{%
   \cite[p.~140]{OP}}]\label{le:0-1}
   For any given $s,t>0$, the following 
   $\sigma$-algebra is trivial:
   \begin{equation}
      \mathfrak{G}(s,t) :=
      \bigcap_{\varepsilon \in\mathbb{Q}_+}
        \sigma \Big\{ W(u,v)-W(s,t):\
        \left| (s,t)-(u,v) \right| <\varepsilon \Big\},
   \end{equation}
   where $|\cdots|$ denotes the $\ell^\infty$-norm on
   $\mathbb{R}^2$ for the sake of concreteness.
\end{lemma}

\begin{proof} \textbf{(Sketch)}
   This requires ideas that are very close to those
   introduced in the proof of
   Corollary~\ref{co:0-1}. Thus, we will 
   indicate only the essential differences between
   the two proofs.
   
   We can think of Brownian sheet as
   the distribution function of white noise.
   Namely, let $\dot{\mathrm{W}}$ 
   denote one-dimensional white noise
   spread over $\mathbb{R}^2$ and \emph{define}
   the Brownian sheet $W$ as
   \begin{equation}\label{white-noise}
      W(u,v) :=\dot{\mathrm{W}}
      \Big( [0,u]\times[0,v] \Big),
      \qquad\forall u,v\ge 0.
   \end{equation}
   [The process $\dot{\mathrm{W}}$ is a
   well-defined vector-valued random measure with
   values in $L^2(\IP)$;
   cf.~\cite[p. 283--285]{walsh:84} and%
   ~\cite[Chapter 5, \S1.3]{MPP} for more details.]
   In this way, we can write
   \begin{equation}
      \mathfrak{G}(s,t) = \bigcap_{\varepsilon\in
      \mathbb{Q}_+} \sigma \left\{
      \dot{\mathrm{W}}\Big(
      [0,s+\varepsilon] \times
      [0,t+\varepsilon] \setminus [0,s-\varepsilon] \times
      [0,t-\varepsilon]\Big) \right\}.
   \end{equation}
   Now suppose that $R$ is
   a rectangle with sides parallel to the axes,
   and that $R$ does not intersect the
   annulus, 
   $[0,s+\varepsilon] \times
   [0,t+\varepsilon]\setminus[0,s-\varepsilon] \times
   [0,t-\varepsilon]$. Then the elementary properties
   of white noise show us that 
   $\dot{\mathrm{W}}(R)$
   is independent of $\mathfrak{G}(s,t)$.
   To finish, consider a finite union of such $R$'s, and
   ``take limits'' in  a manner similar to what we
   did in the proof of Corollary~\ref{co:0-1}.
\end{proof}

\begin{proof}[Proof of Theorem~\ref{th:kendall}]
   We may assume without loss of too much
   generality that $s=t=1$. Consider 
   \begin{equation}
      \mathrm{B}(r):= \left\{ (x,y)\in\mathbb{R}^2:\
      \left| (x,y)-(1,1) \right| \le r \right\},
   \end{equation}
   which is the $\ell^\infty$-ball of radius
   $r>0$ about $(1,1)\in\mathbb{R}^2$.
   Also let $\partial \mathrm{B}(r)$ denote its Euclidean
   boundary; this is the perimeter
   boundary of the square of side $2r$ centered at
   $(1,1)$. Let
   \begin{equation}\label{J}\begin{split}
      \mathsf{J}(r)  := &\left\{\omega:\
         W(1,1) > \sup_{(u,v) \in \partial
         \mathrm{B}(r)} W(u,v) \right\}\\
      =& \left\{\omega:\
         W(1,1) > \sup_{(u,v) \in \partial
         \mathrm{B}(r)} W\left(ru+1-r,rv+1-r\right) \right\}.
   \end{split}\end{equation}
   Then the theorem follows at once from 
   Lemma~\ref{le:0-1} and the following:
   \begin{equation}\label{kendall:goal}
      \liminf_{r\to 0^+} \IP 
      \{ \mathsf{J} (r) \} >0.
   \end{equation}
   For then it follows that with probability one,
   infinitely many of the $\mathsf{J}(n^{-1})$'s 
   occur, and clearly this
   does the job. Therefore, it suffices to prove
   (\ref{kendall:goal}). By
   (\ref{white-noise}), we can write the following
   path decomposition: For any $r\in(0,1)$ fixed,
   \begin{equation}\label{XYZ}\begin{split}
      W(1-r+ur,1-r+vr) = &\sqrt{(1-r)r}\Big[
         X(u)+Y(v)\Big]\\
      &+ rZ(u,v)-W(1-r,1-r),
      \end{split}\quad\forall u,v\ge 0,\qquad
   \end{equation}
   where $X$ and $Y$ are standard Brownian motions,
   $Z$ is a standard Brownian sheet, and the three
   are independent from one another as well as
   from $W(1-r,1-r)$. Indeed, here are the formulas
   for $(X,Y,Z)$ in terms of the white noise $\dot{\mathrm{W}}$
   of (\ref{white-noise}):
   \begin{equation}\label{eq:decompose}\begin{split}
      X(u) := & \frac{1}{\sqrt{(1-r)r}} 
         \dot{\mathrm{W}} \Big( [1-r,1-r+ur]\times[0,1-r]
         \Big)\\
      Y(v) := & \frac{1}{\sqrt{(1-r)r}}
         \dot{\mathrm{W}}\Big( [0,1-r]\times[1-r,1-r+vr]
         \Big)\\
     Z(u,v) := &\frac{1}{r}\dot{\mathrm{W}} \Big(
         [1-r,1-r+ur]\times [1-r,1-r+vr]\Big).
   \end{split}\end{equation}
   We only need these formulas to check the
   assertions about $(X,Y,Z)$, and this
   is only a matter of checking a few covariances. In light of
   the second equality in (\ref{J}), we have
   $\omega\in\mathsf{J}(r)$ if and only if
   \begin{equation}\begin{split}
      &\sqrt{(1-r)r}\Big[ X(1)+Y(1)\Big] + rZ(1,1)\\
      &\hskip1in > \sqrt{(1-r)r} \Big[ X(u)+Y(v) \Big]
         +rZ(u,v),\hskip1in\ 
      \end{split}
   \end{equation}
   simultaneously for all $(u,v)\in\partial \mathrm{B}(1)$.
   But thanks to the sample function continuity
   of $Z$,
   the supremum of $|Z(u,v)|$ over
   all $(u,v)\in\partial \mathrm{B}(1)$ is bounded almost surely.
   Thus, we can divide the preceding display
   by $\sqrt{r}$ and let $r\to 0^+$ to see that
   \begin{equation}
      \lim_{r\to 0^+}\IP\{\mathsf{J}(r)\}
      =\IP\Big\{  \forall (u,v)\in\partial \mathrm{B}(1):\
      X(1)+Y(1) > X(u)+Y(v) \Big\},
   \end{equation}
   and it is easy to see that the latter probability
   is strictly positive. To see this, let
   $X_i:= X(i)$, $Y_i:= Y(i)$ ($i=1,2$). Also let
   $X_2^+:= \sup_{0\le u\le 2}X(u)$, and
   $Y_2^+:= \sup_{0\le v\le 2}Y(v)$, and note that
   the latter probability is equal to
   \begin{equation}\begin{split}
      &\IP \Big\{ X_1 +Y_1 >X_2 + Y_2^+, 
         X_1 + Y_1 >X_2^+ + Y_2 \Big\} \\
      &\ge \IP \Big\{ X_1-X_2 \ge 1, X_1 - X_2^+ \ge -1 \Big\}
         \times \IP \Big\{
         Y_2^+-Y_1 \le 1 , Y_2-Y_1\le -1 \Big\},
   \end{split}\end{equation}
   and this is easily seen to be positive.
   This verifies (\ref{kendall:goal}) and the result
   follows.
\end{proof}

\section{Criterion for Hitting Points}
Thus far, we have only touched upon results that
hold (or do not hold) for quasi-every 
one-dimensional function. This in turn has led us
to the one-dimensional Brownian sheet. Now we turn to
results in higher dimensions. 
With this in mind, let $W$
denote the Brownian sheet in $d$ dimensions
and $\mathfrak{m}_d$ the standard Lebesgue
measure on the Lebesgue-measurable subsets
of $\mathbb{R}^d$. Stated in terms of Wiener measure,
this yields the following.

We begin with the classical fact that
$d$-dimensional Brownian can hit points if
and only if $d=1$.

\begin{theorem}[P.\@ L\'evy]
   The following are equivalent: Given any $x\in\mathbb{R}^d$,
   \begin{enumerate}
      \item[(i)] Almost every continuous
                 $f:[0,1]\to\mathbb{R}^d$
                 avoids $\{x\}$; i.e., $x\not\in f([0,1])$.
      \item[(ii)] Almost every continuous
                 $f:[0,1]\to\mathbb{R}^d$
                 has a Lebesgue-null range; i.e., 
                 $\mathfrak{m}_d(f([0,1]))=0$.
      \item[(iii)] $d\ge 2$.
   \end{enumerate}
\end{theorem}

For the above conditions (i) and (ii) to hold for
quasi-every function $f$, one needs the stronger
condition that $d\ge 4$. Indeed, S.\@ Orey and W.\@ E.\@ Pruitt
have proven the following result.

\begin{theorem}[S.\@ Orey and W.\@ E.\@ Pruitt~
    {\cite[Theorems 3.3 and 3.4]{OP}}]
    \label{thm:hitpoints}
    The following are equivalent: Given any 
    fixed $x\in\mathbb{R}^d$:
    \begin{enumerate}
        \item[(i)] With probability one, $d$-dimensional
             Brownian sheet does not hit $\{x\}$.
        \item[(ii)] The random set $W(\mathbb{R}^2_+)$ 
             has zero $d$-dimensional Lebesgue measure.
        \item[(iii)] $d \ge 4$.
        \item[(iv)] Quasi-every $d$-dimensional
             continuous function avoids $x$. 
    \end{enumerate}
\end{theorem}

Let $U(s,t):= e^{-s/2} W(e^s,t)$ as before, and
note that (iv) is equivalent to the
following:
\begin{equation}
   \IP\left\{ \exists s\in[0,1]:\ \text{ for some }
   t>0,\ U(s,t)=x  \right\}=0,\qquad
   \forall x\in\mathbb{R}^d.
\end{equation}
On the other hand, by the Cameron--Martin 
formula, the law of $\{ U(s,t);\,
s,t\in[0,1]\}$ is mutually absolutely continuous 
with respect to the law of 
$\{ W(s,t);\, s\in[1,e], t\in[0,1]\}$; cf.~\cite{cameron-martin}.
Therefore, (iv)$\Leftrightarrow$(i), and we only
need to prove the equivalence of (i)--(iii).
I will describe most of this proof in three steps.

\begin{proof}[Proof of (i){\boldmath
    $\Rightarrow$}(ii)]
    By Fubini's theorem,
    \begin{equation}
       \IE\left\{ \mathfrak{m}_d \Big( W([0,1]^2) \Big)
       \right\} = \int_{\mathbb{R}^d}
       \IP\left\{ x\in W([0,1]^2) \right\}\, dx=0,
    \end{equation}
    thanks to (i). Scaling then shows
    that with probability one,
    $\mathfrak{m}_d(W(\mathbb{R}^2_+))=0$.
\end{proof}

\begin{proof}[Proof of (ii){\boldmath
    $\Rightarrow$}(iii)]
    We will use the Fourier-analytical ideas 
    of~\cite{kahanebook}, and prove that if (iii)
    fails, then so will (ii).
    Thus, let us assume that $d\le 3$, and consider
    the \emph{occupation (or sojourn) measure},
    \begin{equation}\label{eq:OccMeas}
       \sigma(A) :=
       \int_0^\infty\int_0^\infty
       e^{-s-t} \mathbf{1}_A\left(
       W(s,t) \right) \, ds\, dt.
    \end{equation}
    Its Fourier transform is given by
    \begin{equation}
       \widehat{\sigma}(\xi) = 
       \int_0^\infty\int_0^\infty
       e^{-s-t} e^{i\xi\cdot W(s,t)}
       \, ds\, dt,\qquad\forall\xi\in\mathbb{R}^d.
    \end{equation}
    Our strategy is to show that with probability one,
    $\widehat{\sigma}\in L^2(\mathbb{R}^d)$.
    If so, then by the Plancherel theorem,
    $\sigma$ is a.s. absolutely continuous
    with respect to $\mathfrak{m}_d$, and
    $\left(\frac{d\sigma}{d\mathfrak{m}_d}\right)\in
    L^2(\mathbb{R}^d)$ almost surely. But the fact that
    $\sigma\left(W(\mathbb{R}^2_+) \right)
    =1$ implies that
    \begin{equation}
       \mathfrak{m}_d(W(\mathbb{R}^2_+)) =
       \int_{W(\mathbb{R}^2_+)} \frac{d\sigma (\xi)}{%
       d\mathfrak{m}_d}\, \mathfrak{m}_d(d\xi)
       =1.
    \end{equation}
    Thus,
    (ii)$\Rightarrow$(iii) follows once we
    show that
    $\IE\{\|\widehat{\sigma}
    \|^2_{L^2(\mathbb{R}^d)}\}<+\infty$. The latter
    expectation is equal to the following:
    \begin{equation}\begin{split}
       &\IE \left\{ \| \widehat{\sigma}
          \|^2_{L^2(\mathbb{R}^d)} \right\}\\
       &\quad= \int_{\mathbb{R}^2_+}
          \int_{\mathbb{R}^2_+}
          \int_{\mathbb{R}^d}
          e^{-s_1-s_2-t_1-t_2}
          \IE\left\{ e^{
          i\xi\cdot\left[ W(s_1,s_2) - W(t_1,t_2)
          \right] } \right\}\, d\xi\, ds\, dt\\
       &\quad= \int_{\mathbb{R}^2_+}
          \int_{\mathbb{R}^2_+}
          \int_{\mathbb{R}^d}
          e^{-s_1-s_2-t_1-t_2}
          \exp \left( -\frac{|\xi|^2 \tau^2(s,t)}{2} 
          \right) \, d\xi\, ds\, dt,\hskip.8in
    \end{split}\end{equation}
    where $\tau^2(s,t):= 
    \text{Var}[ W(s_1,s_2) -
    W(t_1,t_2) ]$ is the $\mathfrak{m}_d$-measure
    of the set difference between the rectangles
    $[0,s_1]\times [0,s_2]$ and
    $[0,t_1]\times [0,t_2]$. A picture will convince
    you that no matter how the two said rectangles
    are situated, we always have the bound,
    $\tau^2(s,t) \ge (s_2\wedge t_2) |s_1-t_1|
    +(s_1\wedge t_1) |s_2-t_2|.$
    Thus,
    \begin{equation}\begin{split}
       &\IE \left\{ \| \widehat{\sigma}
          \|^2_{L^2(\mathbb{R}^d)} \right\}\\
       &\displaystyle
          \quad\le 4\iint_{0\le s_1\le t_1}
          \iint_{0\le s_2\le t_2}
          \int_{\mathbb{R}^d}
          e^{-s_1-s_2-t_1-t_2}\\
       &\displaystyle\qquad\times
          \exp\left(-\frac{|\xi|^2 \left[
          s_2(t_1-s_1) + s_1(t_2-s_2)\right]
          }{2} \right)\, d\xi\, ds\, dt\\
       &\displaystyle
          \quad= 4\int_{\mathbb{R}^2_+}
          \int_{\mathbb{R}^2_+}
          \int_{\mathbb{R}^d}
          e^{-2s_1-2s_2-t_1-t_2}
          \exp\left(-\frac{|\xi|^2 \left[
          s_2t_1 + s_1t_2\right]
          }{2} \right)\, d\xi\, ds\, dt.\qquad
    \end{split}\end{equation}
    We can integrate $[dt]$ to deduce that
    $\IE \{ \| \widehat{\sigma}
    \|^2_{L^2(\mathbb{R}^d)} \}\le 4
    \int_{\mathbb{R}^d}
    \mathcal{Q}(\xi) \, d\xi,$
    where
    \begin{equation}
        \mathcal{Q}(\xi) := \int_{\mathbb{R}^2_+}
        e^{-2(s_1+s_2)}\left(1+
        \frac{|\xi|^2s_1}{2}\right)^{-1}\cdot
        \left(1+\frac{|\xi|^2s_2}{2}\right)^{-1}\, ds.
    \end{equation}
    Evidently, $\mathcal{Q}(\xi)$ is bounded, and is
    $O(|\xi|^{-4})$ as $|\xi|\to\infty$. Therefore,
    whenever $d< 4$, then 
    $\mathcal{Q}\in L^1(\mathbb{R}^d)$, and so $\IE\{ \|
    \widehat{\sigma}
    \|^2_{L^2(\mathbb{R}^d)} \}<+\infty$, as asserted.
\end{proof}

\begin{proof}[Partial Proof of (iii){\boldmath$\Rightarrow$}(i)]
    This is the most interesting, as well as difficult, part
    of Theorem~\ref{thm:hitpoints}, and 
    I will present a proof that is valid in the ``supercritical''
    regime $d\ge 5$. When $d=4$, the known proofs are much longer
    and not included here.
    Clearly, it suffices to show that whenever $0<a<b$,
    $\mathrm{P}\{x\in W([a,b]^2)\}=0$. Without loss of too much
    generality, we will do this for $a=1$ and $b=2$.

    Let us fix $\varepsilon\in(0,1)$ and an integer $n\ge 1$,
    and consider the covering $[1,2]^2=\cup_{i,j=0}^n I_{i,j}$,
    where $I_{i,j}:= [1+(i/n),1+(i+1)/n]
    \times [1+(j/n),1+(j+1)/n]$. Now if there exists $(s,t)\in I_{i,j}$
    such that $|W(s,t)-x|\le\varepsilon$, then
    \begin{equation}\begin{split}
       &\left| W\left( 1+\frac{i}{n},1+\frac{j}{n}\right)
          - x\right|\\
       &\quad \le  \varepsilon + \sup_{0\le u,v\le \frac1n}
          \left| W\left( 1+\frac{i}{n},1+\frac{j}{n}\right)
          - W\left( 1+\frac{i}{n}+u,1+\frac{j}{n}+v\right)\right|\\
       &\quad := \varepsilon + \delta_{i,j;n}.
    \end{split}\end{equation}
    Because of the white noise representation~(\ref{white-noise}) of $W$,
    $\delta_{i,j;n}$ is independent of $W(1+in^{-1},1+jn^{-1})$.
    Moreover, the probability density of
    $W(1+in^{-1},1+jn^{-1})$ is uniformly
    bounded above by one. Therefore,
    \begin{equation}\label{eq:supercrit}\begin{split}
       \mathrm{P}\left\{ \exists (s,t)\in I_{i,j}:\
          \left| W(s,t)-x\right| \le \varepsilon \right\}
          & \le C_d\mathrm{E} \left[ \left( \varepsilon + 
          \delta_{i,j;n}\right)^d \right]\\
       &\le 2^d C_d
          \mathrm{E}\left[\varepsilon^d +\delta_{i,j;n}^d \right],
    \end{split}\end{equation}
    where $C_d$ denotes the volume of the unit ball in $\mathbb{R}^d$.
    On the other hand, by the white noise representation of $W$,
    \begin{equation}\begin{split}
       &W\left( 1+\frac{i}{n}+u,1+\frac{j}{n}+v\right)-
          W\left( 1+\frac{i}{n},1+\frac{j}{n}\right)\\
      &\quad =\sqrt{1+\frac{j}{n}}\ B(u) +
          \sqrt{1+\frac{i}{n}}\ B'(v) +Z(u,v),
    \end{split}\end{equation}
    where $B$, $B'$, and $Z$ are independent, $B$ and $B'$
    are Brownian motions, and $Z$ is a Brownian sheet.
    Consequently, we can take absolute values and maximize over
    $u,v\le n^{-1}$ to see that given $0\le i,j\le n^{-1}$,
    \begin{equation}
       \delta_{i,j;n}\le \sqrt2 \sup_{u\in[0,1/n]}
       |B(u)| + \sqrt2 \sup_{v\in[0,1/n]} |B'(v)|
       +\sup_{u,v\in[0,1/n]} |Z(u,v)|.
    \end{equation}
    This and scaling show the existence of a constant
    $K_d$ such that $\mathrm{E}[\delta_{i,j;n}^d]\le K_d n^{-d/2}.$
    Thus, according to~(\ref{eq:supercrit}),
    \begin{equation}
       \mathrm{P}\left\{ \exists (s,t)\in I_{i,j}:\
       \left| W(s,t)-x\right| \le \varepsilon \right\}
       \le  2^d C_d K_d \left[\varepsilon^d + n^{-d/2} \right].
    \end{equation}
    We can sum this over all $0\le i,j\le n$ to see that
    \begin{equation}
       \mathrm{P}\left\{ \exists (s,t)\in [1,2]^2:\
       \left| W(s,t)-x\right| \le \varepsilon \right\}
       \le  2^d C_d K_d (n+1)^2 
       \left[\varepsilon^d + n^{-d/2} \right].
    \end{equation}
    Because this is valid for all $n\ge 1$, we can choose $n := \lfloor
    \varepsilon^{-2}\rfloor$, and deduce that
    \begin{equation}
       \sup_{x\in\mathbb{R}^d}\mathrm{P}\left\{ \exists (s,t)\in [1,2]^2:\
       \left| W(s,t)-x\right| \le \varepsilon \right\}
       = O\left( \varepsilon^{d-4} \right),\qquad(\varepsilon\to 0).
    \end{equation}
    In particular, if $d\ge 5$, then with 
    probability one, $x\not\in W([1,2]^2)$, as claimed.
\end{proof}

\section{The O-U Process on Wiener Space and Two Open Problems}

In this section, I very briefly sketch the connection between
the process $Y$ and symmetric forms. This is the starting point
for the introduction of the methods of potential theory; an 
area which is not the focus of the present article. I then conclude
this discussion by presenting two of my favorite 
open problems in this general area.

We have seen already that the Ornstein--Uhlenbeck process
$Y$ of (\ref{Y}) is a stationary diffusion on $\Omega$
(the space of all real continuous functions on $[0,1]$)
whose stationary measure is Wiener's measure. Moreover, in
the usual notation of Markov processs, we have
the following for all continuous $\phi:\Omega\to\mathbb{R}_+$,
$s>0$, and $x\in\Omega$:
\begin{equation}\begin{split}
   T_s \phi (x) & :=  \IE
      \big[ \phi(Y_s)\, \big|\, Y_0=x\big] =
      \IE \left[ \phi \Big( e^{-s/2} W(e^s,\bullet) \Big)\,\Big|\,
      W(1,\bullet)=x\right]\\
   &=\IE \left[ \phi \Big( e^{-s/2} \Big[
      W(e^s,\bullet) - W(1,\bullet)\Big] +e^{-s/2}x\Big)\right]\\
   &=\IE \left[ \phi \Big( e^{-s/2} \Big[
      W(e^s-1,\bullet) \Big] +e^{-s/2}x\Big)\right]\\
   &=\IE \left[ \phi \Big( \sqrt{1-e^{-s}}
      W(1,\bullet) +e^{-s/2}x\Big)\right].
\end{split}\end{equation}
On the other hand, $W(1,\bullet)$ is just a Brownian motion,
and this leads to the following.

\begin{lemma}[Mehler's Formula]
   If $\mu$ denote the Wiener measure on 
   the classical Wiener space $(\Omega,\mathcal{B}(\Omega))$,
   then the transition semigroup of the diffusion $Y$
   are given by the operator-formula: For all $s>0$ and all
   continuous functions $x\in\Omega$,
   \begin{equation}
      T_s f(x) = \int_\Omega f\left(
      \sqrt{1-e^{-s}} y + e^{-s/2}x \right)\, \mu(dy).
   \end{equation}
\end{lemma}It is also a simple matter to check
that $T_s$ is a symmetric semigroup on
the classical Wiener space; i.e., that 
$\left\langle g,T_s f\right\rangle_\Omega
=\left\langle T_s g,f\right\rangle_\Omega$,
where $\langle u,v\rangle_\Omega$ is the 
covariance form, $\int_\Omega uv\, d\mu$.

Thus, the standard theory of symmetric Markov processes
constructs a Dirichlet form $\mathcal{E}$ for $Y$
killed at an exponential rate that is formally
defined as follows:
\begin{equation}
   \mathcal{E}(f,g) := \lim_{s\to 0^+}
   \left\langle  s^{-1} ( f -T_s f )
   ,g\right\rangle_\Omega +\langle
   f ,g\rangle_\Omega;
\end{equation}
cf. \cite[(1.3.15) and Theorem 1.4]{fukushima-et-al}.
Now given any open set $G\subseteq\Omega$
one has the following identity; it relates the capacity
of \S1 to the Dirichlet forms of \cite{fukushima-et-al}:
\begin{equation}\label{cap}
   \text{\rm Cap}(G) = \inf \Big\{
   \mathcal{E}(f,f);\, f\in
   \text{\rm dom}(\mathcal{E}),\
   f\ge 1\ \mu\text{-a.e. on }G
   \Big\},
\end{equation}
where $\text{\rm dom}(\mathcal{E})$ denotes the domain
of $\mathcal{E}$; i.e., all $f\in L^2(\Omega)$
such that $\mathcal{E}(f,f)<+\infty.$
Furthermore, for a general set $A\subseteq\Omega$,
$\text{\rm Cap}(A)=\inf\{\text{\rm Cap}(G):\
A\subseteq G\text{ open},\}$.
The latter remarks are proved in 
\cite[p. 164]{fukushima},
and are another starting point for the analytic
treatment of many of the 
quasi-sure results within the references.

As promised earlier, we conclude this paper by
presenting two open problems:
\begin{enumerate}
   \item[{\bf OP-1}] 
      P.\@ Malliavin has introduced a parametric
      family of Gaussian capacities, one of which is the
      $\text{\rm Cap}$ of this and the first section. 
      Is there a ``truly probabilistic'' description
      (i.e., one involving concrete random processes)
      of all of these capacities? If so, do the quasi-sure
      results of this paper continue to holds if $\text{\rm
      Cap}$ is everywhere replaced by any and all
      of the said capacities? For related results, see
      the results of M.\@ Takeda (\cite{takeda}).
   \item[{\bf OP-2}]
      One of the outstanding open problems of the geometry
      of Brownian sheet is the following: Let $\mathcal{L}$
      denote the complement of the 
      zero-set of an $N$-parameter Brownian sheet. Does
      the complement of $\mathcal{L}$ have an infinite
      connected component? When $N=2$, this was answered
      in the negative
      by W.\@ S.\@ Kendall (\cite{kendall}), but the general question
      remains open.
\end{enumerate}

\end{document}